\documentclass[11pt,a4paper]{article}
	\usepackage{amsmath,amssymb}
	\newtheorem{lemma}{Lemma}
	\newtheorem{example}{Example}
	\newtheorem{definition}{Definition}
	\newtheorem{theorem}{Theorem}
	\newtheorem{corollary}{Corollary}
	\newtheorem{conjecture}{Conjecture}
	\newcommand{\eop}{\hfill{$\Box$}}
	\newenvironment{proof}
	{\begin{trivlist}\item[]{{\sc Proof.}}}{\eop\noindent\end{trivlist}}

\begin{document}
	\title{Convex hulls of polyominoes}
	\author{{\sc Sascha Kurz}\thanks{sascha.kurz@uni-bayreuth.de, http://www.wm.uni-bayreuth.de/index.php?id=sascha}\\ 
		  Business Mathematics, University of Bayreuth\\ 
		  D-95440 Bayreuth, Germany}
	\maketitle		 
	\vspace*{-4mm}
	\noindent
	{
		\center\small{Keywords: polyominoes, convex hull, dido-type problem, isoperimetric inequality\\ 
		 MSC: 05B50$^\star$, 05D99, 52C99\\}
	}
	\noindent
	\rule{\textwidth}{0.3 mm}

	\begin{abstract}
		\noindent
		In this article we prove a conjecture of Bezdek, Bra\ss, and Harborth concerning the maximum volume of the convex hull 
		of any facet-to-facet connected system of $n$ unit hypercubes in $\mathbb{R}^d$ \cite{Bezdek1994}. For $d=2$ we enumerate 
		the extremal polyominoes and determine the set of possible areas of the convex hull for each $n$.
	\end{abstract}
	\noindent
	\rule{\textwidth}{0.3 mm}

\section{Introduction}

In the legend \cite{legend} of the founding of Carthage, Queen Dido purchased the right to
get as much land as she could enclose with the skin of an ox. She splitted the
skin into thin stripes and tied them together. Using the natural boundary of the
sea and by constructing a giant semicircle she enclosed more land than the seller
could have ever imagined.

Dido-type problems have been treated by many authors i.e. \cite{0703.05036,Bezdek1994,0159.24203,
dido_5,0999.52007}, here we consider the maximum 
volume of a union of unit hypercubes. A $d$-dimensional polyomino is a facet-to-facet
connected system of $d$-dimensional	unit hypercubes. Examples for $2$-dimensional 
polyominoes are the pieces of the computer game Tetris.

In 1994 Bezdek, Bra\ss, and Harborth conjectured that the maximum volume of the convex
hull of a $d$-dimensional polyomino consisting of $n$ hypercubes is given by
$$
	\sum\limits_{I\subseteq\{1,\dots,d\}}\frac{1}{|I|!}\prod\limits_{i\in
	I}\left\lfloor\frac{n-2+i}{d}\right\rfloor,
$$
but were only able to prove it for $d=2$. In Section \ref{sec_3} we prove this conjecture. 
They also asked for the number $c_2(n)$ of different polyominoes with $n$ cells and maximum area 
$n+\left\lfloor\frac{n-1}{2}\right\rfloor\left\lfloor\frac{n}{2}\right\rfloor$. 
In Section \ref{sec_2} we prove

\begin{theorem}
	\label{thm_count_2}
	$$
		c_2(n)=\left\{
			\begin{array}{rcl}
				\frac{n^3-2n^2+4n}{16} & \text{if} & n\equiv 0\mod 4,\\
				\frac{n^3-2n^2+13n+20}{32} & \text{if} & n\equiv 1\mod 4,\\
				\frac{n^3-2n^2+4n+8}{16} & \text{if} & n\equiv 2\mod 4,\\
				\frac{n^3-2n^2+5n+8}{32} & \text{if} & n\equiv 3\mod 4.\\
			\end{array}
		\right.
	$$
\end{theorem}

Besides the maximum area $n+\left\lfloor\frac{n-1}{2}\right\rfloor\left\lfloor\frac{n}{2}\right\rfloor$ and 
the minimum area $n$ of the convex hull of polyominoes with $n$ cells several other values may be attained. For 
each $n$ we characterize the corresponding sets.

\begin{theorem}
	\label{thm_possible_values}
	A polyomino consisting of $n$ cells with area $\alpha=n+\frac{m}{2}$ of the convex hull exists if and 
	only if $m\in\mathbb{N}_0$, $0\le m\le \left\lfloor\frac{n-1}{2}\right\rfloor\left\lfloor\frac{n}{2}\right\rfloor$, 
	and $m\neq 1$ if $n+1$ is a prime.
\end{theorem}

\section{The planar case}
\label{sec_2}

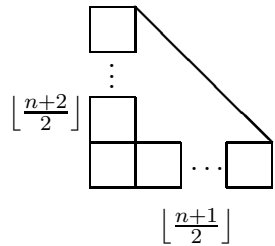
\begin{figure}[h]
	\begin{center}
		\setlength{\unitlength}{0.6cm}
		\begin{picture}(5,5)
			\put(1,1){\line(1,0){2}}
			\put(1,1){\line(0,1){2}}
			\put(2,2){\line(0,-1){1}}
			\put(2,2){\line(-1,0){1}}
			\put(2,3){\line(-1,0){1}}
			\put(2,3){\line(0,-1){1}}
			\put(3,2){\line(-1,0){1}}
			\put(3,2){\line(0,-1){1}}
			\put(4,1){\line(1,0){1}}
			\put(4,1){\line(0,1){1}}
			\put(5,2){\line(0,-1){1}}
			\put(5,2){\line(-1,0){1}}
			\put(1,4){\line(1,0){1}}
			\put(1,4){\line(0,1){1}}
			\put(2,5){\line(0,-1){1}}
			\put(2,5){\line(-1,0){1}}
			\put(3.2,1.4){$\dots$}
			\put(1.4,3.2){$\vdots$}
			\thicklines
			\put(2,5){\line(1,-1){3}}
			\put(-0.8,2.5){$\left\lfloor\frac{n+2}{2}\right\rfloor$}
			\put(2.5,0){$\left\lfloor\frac{n+1}{2}\right\rfloor$}
		\end{picture}
		\caption{2-dimensional polyomino with maximum convex hull.}
		\label{fig_extremal_2}
	\end{center}
\end{figure}

An example which attains the upper bound $n+\left\lfloor\frac{n-1}{2}\right\rfloor\left\lfloor\frac{n}{2}\right\rfloor$ for 
the area of the convex hull of a polyomino with $n$ cells is quite obvious, see Figure \ref{fig_extremal_2}. Instead 
of proving this upper bound by induction over $n$ we specify polyominoes by further parameters and then apply an induction argument.

We describe these parameters for the more general $d$-dimensional case and therefore denote the standard 
coordinate axes of $\mathbb{R}^d$ by $1,\dots,d$. Every $d$-dimensional polyomino has a smallest surrounding box 
with side lengths $l_1,\dots,l_d$, where $l_i$ is the length in direction $i$. If we build up a polyomino cell 
by cell then after adding a cell one of the $l_i$ will increase by $1$ or none of the $l_i$ will increase. In the 
second case we increase $v_i$ by $1$, where the new hypercube has a facet-neighbor in direction of axis $i$. If $M$ 
is the set of axis-directions of facet-neighbors of the new hypercube, then we will increase $v_i$ by $1$ for only 
one $i\in M$. Since at this position there is the possibility to choose, we must face the fact that there might 
be different tuples $(l_1,\dots,l_d,v_1,\dots,v_d)$ for the same polyomino. We define $v_1=\dots=v_d=0$ for the polyomino 
consisting of a single hypercube. This definition of the $l_i$ and the $v_i$ leads to 
\begin{equation}
	\label{main_equation}
	n = 1 + \sum_{i=1}^d (l_i-1) + \sum_{i=1}^d v_i .
\end{equation}

\begin{example}
	The possible tuples describing a rectangular $2\times 3$-polyomino are $(2,3,2,0)$, $(2,3,1,1)$, and $(2,3,0,2)$.
\end{example}

\begin{definition}
	\begin{eqnarray*}
		f_2(l_1,l_2,v_1,v_2) & = &1+(l_1-1)+(l_2-1)+\frac{(l_1-1)(l_2-1)}{2}\\
		&&+v_1+v_2+\frac{v_1(l_2-1)}{2}+\frac{v_2(l_1-1)}{2}+\frac{v_1v_2}{2}.
	\end{eqnarray*}
\end{definition}

\begin{lemma}
	\label{lemma_two_dimensional}
	The area of the convex hull of a $2$-dimensional polyomino with tuple $(l_1,l_2,v_1,v_2)$ is at most $f_2(l_1,l_2,v_1,v_2)$.
\end{lemma}
\begin{proof}
	We prove the statement by induction on $n$, using equation \ref{main_equation}. For $n=1$ only $l_1=l_2=1$, $v_1=v_2=0$ is 
	possible. With $f_2(1,1,0,0)=1$ the induction base is done. Now we assume that the statement is true for all possible 
	tuples $(l_1,l_2,v_1,v_2)$ with $1 + \sum_{i=1}^d (l_i-1) + \sum_{i=1}^d v_i=n-1$. 
	
	\medskip
	
	Due to symmetry we consider only the growth of $l_1$ or $v_1$, and the area $a$ of the convex hull by adding the $n$-th square.
	
	\begin{enumerate}
		\item[(i)] $l_1$ increases by one:
				
				\begin{minipage}[t]{4.1 cm}
					\setlength{\unitlength}{0.38cm}
					\begin{picture}(10,11)
						\put(1,8){\line(1,0){4}}
						\put(1,10){\line(1,0){4}}
						\put(1,8){\line(0,1){2}}
						\put(3,8){\line(0,1){2}}
						\put(5,8){\line(0,1){2}}
						\put(3,8){\line(0,-1){6}}
						\thicklines
						\put(5,8){\line(-1,-3){2}}
						\thinlines
						\put(3,10){\line(-1,1){4}}
						\thicklines
						\put(5,10){\line(-3,2){6}}
						\thinlines
						\put(3,8){\line(1,1){2}}
						\put(3,9){\line(1,1){1}}
						\put(4,8){\line(1,1){1}}
						\put(6.6,2){\line(0,1){12}}
						\put(6.6,2){\line(-1,0){0.6}}
						\put(6.6,14){\line(-1,0){0.6}}
						\put(7,8){$\le l_2$}
						\stepcounter{figure}
						\put(-1,0){Figure 2: Increasing $l_1$.}
					\end{picture}
				\end{minipage}
				\begin{minipage}[b]{7.52 cm}
					We depict (see Figure 2) the new square by $3$ diagonal lines. Since $l_1$ increases 
					the new square must have a left or a right neighbor. Without loss of generality it has 
					a left neighbor. The new square contributes at most $2$ (thick) lines to the convex hull 
					of the polyomino. By drawing lines from the neighbor square to the endpoints of the new 
					lines we see that the growth is at most $1+\frac{l_2-1}{2}$, a growth of $1$ for the new 
					square and the rest for the triangles. Since 
					$f_2(l_1+1,l_2,v_1,v_2)-f_2(l_1,l_2,v_1,v_2)=1+\frac{l_2-1}{2}+\frac{v_2}{2}$
    				the induction step follows.
  				\end{minipage}
				\newpage
		\item[(ii)] $v_1$ increases by one:
				
				\begin{minipage}[t]{5.1 cm}
					\setlength{\unitlength}{0.38cm}
					\begin{picture}(14,14)
						\put(1.5,8){\line(1,0){4}}
						\put(1.5,10){\line(1,0){4}}
						\put(1.5,8){\line(0,1){2}}
						\put(3.5,8){\line(0,1){2}}
						\put(5.5,8){\line(0,1){2}}
						\put(3.5,10){\line(1,-1){6}}
						\thicklines
						\put(5.5,10){\line(2,-3){4}}
						\thinlines
						\put(3.5,10){\line(-1,1){4}}
						\thicklines
						\put(5.5,10){\line(-3,2){6}}
						\thinlines
						\put(3.5,8){\line(1,1){2}}
						\put(3.5,9){\line(1,1){1}}
						\put(4.5,8){\line(1,1){1}}
						\put(10.6,2){\line(0,1){12}}
						\put(10.6,2){\line(-1,0){0.6}}
						\put(10.6,14){\line(-1,0){0.6}}
						\put(11,8){$\le l_2$}
						\put(3.5,2){\line(1,0){2}}
						\put(3.5,2){\line(0,1){2}}
						\put(5.5,4){\line(-1,0){2}}
						\put(5.5,4){\line(0,-1){2}}
						\put(5.9,3){$\dots$}
						\put(7.5,2){\line(1,0){2}}
						\put(7.5,2){\line(0,1){2}}
						\put(9.5,4){\line(-1,0){2}}
						\put(9.5,4){\line(0,-1){2}}
						\stepcounter{figure}
						\put(-0.5,0){Figure 3: Increasing $v_1$.}
					\end{picture}
				\end{minipage}
				\begin{minipage}[b]{6.5 cm}
					In Figure 3 we depict the new square by $3$ diagonal lines. Without loss of generality
					we assume that the new square has a left neighbor, and contributes at most $2$ lines to 
					the convex hull of the polyomino. As $l_1$ is not increased there must be a square in the 
					same column as the new square. Similar to (i) we draw lines from the neighbor square to the 
					endpoints of the new lines and see that the growth of the area of the convex hull
    				is less than $\frac{l_2-1}{2}$. 
				\end{minipage}
	\end{enumerate}
\end{proof}

\begin{theorem}
	\label{thm_maximum_2}
	The area of the convex hull of a $2$-dimensional polyomino with $n$ unit squares is at most 
	$n+\left\lfloor\frac{n-1}{2}\right\rfloor\left\lfloor\frac{n}{2}\right\rfloor$.
\end{theorem}
\begin{proof}
	For given $n$ we determine the maximum of $f_2(l_1,l_2,v_1,v_2)$. Since $f_2(l_1+1,l_2,v_1-1,v_2)-f_2(l_1,l_2,v_1,v_2)=0$ 
	and due to symmetry we assume $v_1=v_2=0$ and $l_1\le l_2$. With
	$$
		f_2(l_1+1,l_2-1,0,0)-f_2(l_1,l_2,0,0)=\frac{l_2-l_1-1}{2}>0
	$$
	we conclude $0\le l_2-l_1\le 1$. Using equation \ref{main_equation} gives $l_1=\left\lfloor\frac{n+1}{2}\right\rfloor$, 
	$l_2=\left\lfloor\frac{n+2}{2}\right\rfloor$. Thus by inserting in Lemma \ref{lemma_two_dimensional} we receive
	$f_2(l_1,l_2,v_1,v_2)\le n+\frac{1}{2}\left\lfloor\frac{n-1}{2}\right\rfloor\left\lfloor\frac{n}{2}\right\rfloor$. This 
	maximum is attained for example by the polyomino in Figure \ref{fig_extremal_2}.
\end{proof}

In the next lemma we describe the shape of the $2$-dimensional polyominoes with maximum area of the convex hull in order to determine their number $c_2(n)$. 

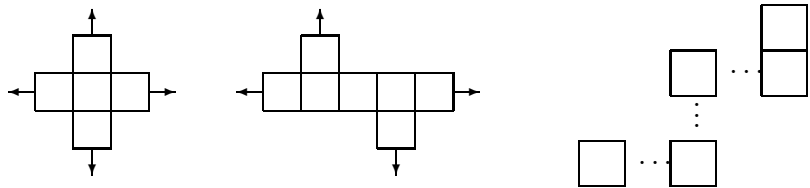
\begin{figure}[h]
	\begin{center}
		\setlength{\unitlength}{5mm}
		\begin{picture}(13,5)
			\put(1,2.5){\vector(-1,0){0.7}}
			\put(1,2){\line(1,0){3}}
			\put(1,3){\line(1,0){3}}
			\put(1,2){\line(0,1){1}}
			\put(2,2){\line(0,1){1}}
			\put(3,2){\line(0,1){1}}
			\put(4,2){\line(0,1){1}}
			\put(4,2.5){\vector(1,0){0.7}}
			\put(2,1){\line(1,0){1}}
			\put(2,1){\line(0,1){1}}
			\put(3,2){\line(0,-1){1}}
			\put(2.5,1){\vector(0,-1){0.7}}
			\put(2,3){\line(0,1){1}}
			\put(3,3){\line(0,1){1}}
			\put(2,4){\line(1,0){1}}
			\put(2.5,4){\vector(0,1){0.7}}
			\put(7,2){\line(1,0){5}}
			\put(7,3){\line(1,0){5}}
			\put(7,2){\line(0,1){1}}
			\put(8,2){\line(0,1){1}}
			\put(9,2){\line(0,1){1}}
			\put(10,2){\line(0,1){1}}
			\put(11,2){\line(0,1){1}}
			\put(12,2){\line(0,1){1}}
			\put(8,3){\line(0,1){1}}
			\put(9,3){\line(0,1){1}}
			\put(8,4){\line(1,0){1}}
			\put(10,2){\line(0,-1){1}}
			\put(11,2){\line(0,-1){1}}
			\put(10,1){\line(1,0){1}}
			\put(7,2.5){\vector(-1,0){0.7}}
			\put(12,2.5){\vector(1,0){0.7}}
			\put(8.5,4){\vector(0,1){0.7}}
			\put(10.5,1){\vector(0,-1){0.7}}
		\end{picture}
		$\quad\quad$
		\setlength{\unitlength}{0.3cm}
		\begin{picture}(10,8)
			\put(0,0){\line(1,0){2}}
			\put(2,0){\line(0,1){2}}
			\put(0,0){\line(0,1){2}}
			\put(0,2){\line(1,0){2}}
			\put(2.6,1){$\dots$}
			\put(4,0){\line(1,0){2}}
			\put(6,0){\line(0,1){2}}
			\put(4,0){\line(0,1){2}}
			\put(4,2){\line(1,0){2}}
			\put(5,2.6){$\vdots$}
			\put(4,4){\line(1,0){2}}
			\put(6,4){\line(0,1){2}}
			\put(4,4){\line(0,1){2}}
			\put(4,6){\line(1,0){2}}
			\put(6.6,5){$\dots$}
			\put(8,4){\line(1,0){2}}
			\put(10,4){\line(0,1){2}}
			\put(8,4){\line(0,1){2}}
			\put(8,6){\line(1,0){2}}
			\put(8,6){\line(1,0){2}}
			\put(10,6){\line(0,1){2}}
			\put(8,6){\line(0,1){2}}
			\put(8,8){\line(1,0){2}}
		\end{picture}
	
		\caption{The two shapes of polyominoes with maximum area of the convex hull and a forbidden sub-polyomino.}
		\label{fig_forbidden}
	\end{center}
\end{figure}	

\begin{lemma}
	\label{lemma_shape_2}
	Every $2$-dimensional polyomino with parameters $l_1$, $l_2$, $v_1$, $v_2$, and with the maximum area 
	$n+\frac{1}{2}\left\lfloor\frac{n-1}{2}\right\rfloor\left\lfloor\frac{n}{2}\right\rfloor$ of the convex hull consists 
	of a linear strip with at most one orthogonal linear strip on each side (see the left two pictures in Figure \ref{fig_forbidden}). 
	Additionally we have $v_1=v_2=0$ and the area of the convex hull is given by $f_2(l_1,l_2,v_1,v_2)$.
\end{lemma}
\begin{proof}
	From the proof of Lemma \ref{lemma_two_dimensional} we deduce $v_1=v_2=0$ and that every sub-polyomino has also the 
	maximum area of the convex hull. Since the area of the polyomino on the right hand side of Figure \ref{fig_forbidden} 
	has an area of the convex hull which is less than $f_2(l_1,l_2,0,0)$ it is a forbidden sub-polyomino and only the 
	described shapes remain. All these polyominoes attain the maximum $f_2(l_1,l_2,0,0)$.
\end{proof}

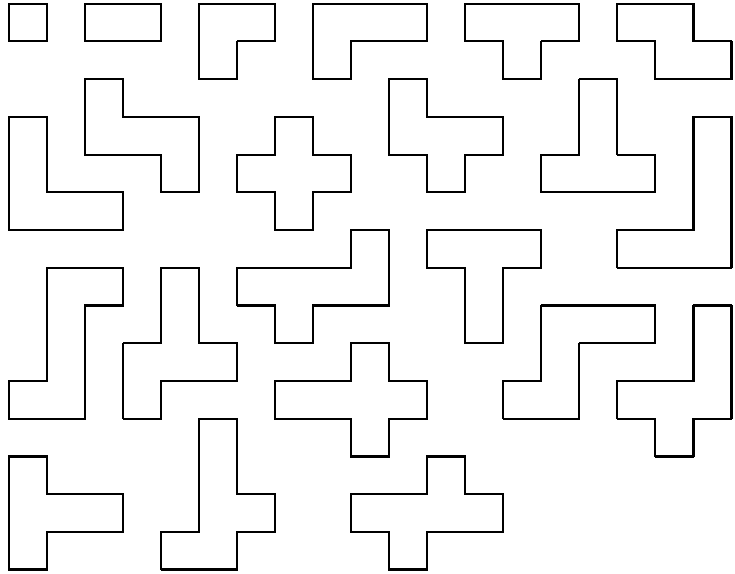
\begin{figure}[h]
	\begin{center}
		\setlength{\unitlength}{0.5cm}
		\begin{picture}(19,15)
			\put(0,0){\line(1,0){1}}
			\put(0,0){\line(0,1){3}}
			\put(0,3){\line(1,0){1}}
			\put(1,3){\line(0,-1){1}}
			\put(1,0){\line(0,1){1}}
			\put(1,1){\line(1,0){2}}
			\put(1,2){\line(1,0){2}}
			\put(3,1){\line(0,1){1}}
			\put(4,0){\line(1,0){2}}
			\put(4,0){\line(0,1){1}}
			\put(6,0){\line(0,1){1}}
			\put(6,1){\line(1,0){1}}
			\put(7,1){\line(0,1){1}}
			\put(7,2){\line(-1,0){1}}
			\put(6,2){\line(0,1){2}}
			\put(6,4){\line(-1,0){1}}
			\put(5,4){\line(0,-1){3}}
			\put(5,1){\line(-1,0){1}}
			\put(10,0){\line(1,0){1}}
			\put(10,0){\line(0,1){1}}
			\put(11,0){\line(0,1){1}}
			\put(11,1){\line(1,0){2}}
			\put(13,1){\line(0,1){1}}
			\put(13,2){\line(-1,0){1}}
			\put(12,2){\line(0,1){1}}
			\put(12,3){\line(-1,0){1}}
			\put(11,3){\line(0,-1){1}}
			\put(11,2){\line(-1,0){2}}
			\put(9,2){\line(0,-1){1}}
			\put(9,1){\line(1,0){1}}
			\put(0,4){\line(1,0){2}}
			\put(0,4){\line(0,1){1}}
			\put(2,4){\line(0,1){3}}
			\put(2,7){\line(1,0){1}}
			\put(3,7){\line(0,1){1}}
			\put(3,8){\line(-1,0){2}}
			\put(1,8){\line(0,-1){3}}
			\put(1,5){\line(-1,0){1}}
			\put(3,9){\line(-1,0){3}}
			\put(3,9){\line(0,1){1}}
			\put(3,10){\line(-1,0){2}}
			\put(1,10){\line(0,1){2}}
			\put(1,12){\line(-1,0){1}}
			\put(0,12){\line(0,-1){3}}
			\put(0,14){\line(1,0){1}}
			\put(0,14){\line(0,1){1}}
			\put(1,15){\line(-1,0){1}}
			\put(1,15){\line(0,-1){1}}
			\put(2,14){\line(1,0){2}}
			\put(2,14){\line(0,1){1}}
			\put(4,15){\line(-1,0){2}}
			\put(4,15){\line(0,-1){1}}
			\put(5,15){\line(1,0){2}}
			\put(5,15){\line(0,-1){2}}
			\put(7,15){\line(0,-1){1}}
			\put(7,14){\line(-1,0){1}}
			\put(6,14){\line(0,-1){1}}
			\put(6,13){\line(-1,0){1}}
			\put(8,15){\line(1,0){3}}
			\put(8,15){\line(0,-1){2}}
			\put(11,15){\line(0,-1){1}}
			\put(11,14){\line(-1,0){2}}
			\put(9,14){\line(0,-1){1}}
			\put(9,13){\line(-1,0){1}}
			\put(12,15){\line(1,0){3}}
			\put(12,15){\line(0,-1){1}}
			\put(15,15){\line(0,-1){1}}
			\put(15,14){\line(-1,0){1}}
			\put(14,14){\line(0,-1){1}}
			\put(14,13){\line(-1,0){1}}
			\put(13,13){\line(0,1){1}}
			\put(13,14){\line(-1,0){1}}
			\put(16,15){\line(1,0){2}}
			\put(18,15){\line(0,-1){1}}
			\put(18,14){\line(1,0){1}}
			\put(19,14){\line(0,-1){1}}
			\put(19,13){\line(-1,0){2}}
			\put(17,13){\line(0,1){1}}
			\put(17,14){\line(-1,0){1}}
			\put(16,14){\line(0,1){1}}
			\put(19,12){\line(0,-1){4}}
			\put(19,8){\line(-1,0){3}}
			\put(16,8){\line(0,1){1}}
			\put(16,9){\line(1,0){2}}
			\put(18,9){\line(0,1){3}}
			\put(18,12){\line(1,0){1}}
			\put(19,7){\line(0,-1){3}}
			\put(19,4){\line(-1,0){1}}
			\put(18,4){\line(0,-1){1}}
			\put(18,3){\line(-1,0){1}}
			\put(17,3){\line(0,1){1}}
			\put(17,4){\line(-1,0){1}}
			\put(16,4){\line(0,1){1}}
			\put(16,5){\line(1,0){2}}
			\put(18,5){\line(0,1){2}}
			\put(18,7){\line(1,0){1}}
			\put(3,4){\line(1,0){1}}
			\put(3,4){\line(0,1){2}}
			\put(4,4){\line(0,1){1}}
			\put(4,5){\line(1,0){2}}
			\put(6,5){\line(0,1){1}}
			\put(6,6){\line(-1,0){1}}
			\put(5,6){\line(0,1){2}}
			\put(5,8){\line(-1,0){1}}
			\put(4,8){\line(0,-1){2}}
			\put(4,6){\line(-1,0){1}}
			\put(7,4){\line(0,1){1}}
			\put(7,4){\line(1,0){2}}
			\put(9,4){\line(0,-1){1}}
			\put(9,3){\line(1,0){1}}
			\put(10,3){\line(0,1){1}}
			\put(10,4){\line(1,0){1}}
			\put(11,4){\line(0,1){1}}
			\put(11,5){\line(-1,0){1}}
			\put(10,5){\line(0,1){1}}
			\put(10,6){\line(-1,0){1}}
			\put(9,6){\line(0,-1){1}}
			\put(9,5){\line(-1,0){2}}
			\put(13,4){\line(1,0){2}}
			\put(15,4){\line(0,1){2}}
			\put(15,6){\line(1,0){2}}
			\put(17,6){\line(0,1){1}}
			\put(17,7){\line(-1,0){3}}
			\put(14,7){\line(0,-1){2}}
			\put(14,5){\line(-1,0){1}}
			\put(13,5){\line(0,-1){1}}
			\put(12,6){\line(1,0){1}}
			\put(13,6){\line(0,1){2}}
			\put(13,8){\line(1,0){1}}
			\put(14,8){\line(0,1){1}}
			\put(14,9){\line(-1,0){3}}
			\put(11,9){\line(0,-1){1}}
			\put(11,8){\line(1,0){1}}
			\put(12,8){\line(0,-1){2}}
			\put(10,9){\line(-1,0){1}}
			\put(9,9){\line(0,-1){1}}
			\put(9,8){\line(-1,0){3}}
			\put(6,8){\line(0,-1){1}}
			\put(6,7){\line(1,0){1}}
			\put(7,7){\line(0,-1){1}}
			\put(7,6){\line(1,0){1}}
			\put(8,6){\line(0,1){1}}
			\put(8,7){\line(1,0){2}}
			\put(10,7){\line(0,1){2}}
			\put(9,10){\line(0,1){1}}
			\put(9,11){\line(-1,0){1}}
			\put(8,11){\line(0,1){1}}
			\put(8,12){\line(-1,0){1}}
			\put(7,12){\line(0,-1){1}}
			\put(7,11){\line(-1,0){1}}
			\put(6,11){\line(0,-1){1}}
			\put(6,10){\line(1,0){1}}
			\put(7,10){\line(0,-1){1}}
			\put(7,9){\line(1,0){1}}
			\put(8,9){\line(0,1){1}}
			\put(8,10){\line(1,0){1}}
			\put(10,11){\line(0,1){2}}
			\put(10,13){\line(1,0){1}}
			\put(11,13){\line(0,-1){1}}
			\put(11,12){\line(1,0){2}}
			\put(13,12){\line(0,-1){1}}
			\put(13,11){\line(-1,0){1}}
			\put(12,11){\line(0,-1){1}}
			\put(12,10){\line(-1,0){1}}
			\put(11,10){\line(0,1){1}}
			\put(11,11){\line(-1,0){1}}
			\put(14,10){\line(1,0){3}}
			\put(17,10){\line(0,1){1}}
			\put(17,11){\line(-1,0){1}}
			\put(16,11){\line(0,1){2}}
			\put(16,13){\line(-1,0){1}}
			\put(15,13){\line(0,-1){2}}
			\put(15,11){\line(-1,0){1}}
			\put(14,11){\line(0,-1){1}}
			\put(4,10){\line(1,0){1}}
			\put(5,10){\line(0,1){2}}
			\put(5,12){\line(-1,0){2}}
			\put(3,12){\line(0,1){1}}
			\put(3,13){\line(-1,0){1}}
			\put(2,13){\line(0,-1){2}}
			\put(2,11){\line(1,0){2}}
			\put(4,11){\line(0,-1){1}}
		\end{picture}
		\caption{Complete set of extremal polyominoes with $n\le 6$ cells.}
		\label{fig_extremal_polyominoes}
	\end{center}
\end{figure}

\textbf{Theorem \ref{thm_count_2}}
$$
	c_2(n)=\left\{
		\begin{array}{rcl}
			\frac{n^3-2n^2+4n}{16} & \text{if} & n\equiv 0\mod 4,\\
			\frac{n^3-2n^2+13n+20}{32} & \text{if} & n\equiv 1\mod 4,\\
			\frac{n^3-2n^2+4n+8}{16} & \text{if} & n\equiv 2\mod 4,\\
			\frac{n^3-2n^2+5n+8}{32} & \text{if} & n\equiv 3\mod 4.\\
		\end{array}
	\right.
$$

\medskip

{\sc Proof of Theorem \ref{thm_count_2}.} (Formula for $c_2(n)$.)
	
We use Lemma \ref{lemma_shape_2} and do a short calculation applying the lemma of Cauchy-Frobenius.
\hfill{$\square$}

\begin{corollary}
	The ordinary generating function for $c_2(n)$ is given by
	$$
		\frac{1+x-x^2-x^3+2x^5+8x^6+2x^7+4x^8+2x^9-x^{10}+x^{12}}{(1-x^2)^2(1-x^4)^2}.
	$$
\end{corollary}

We have depicted the polyominoes with at most $6$ cells and maximum area of the convex hull in Figure \ref{fig_extremal_polyominoes}. For more cells we give only a few concrete numbers:
\begin{eqnarray*}
	\left(c_2(n)\right)_{n=1,\dots} = 1, 1, 1, 3, 5, 11, 9, 26, 22, 53, 36, 93, 64, 151, 94, 228, 143, 329,\\
	195, 455, 271, 611, 351, 798, 460, 1021, 574, 1281, 722, 1583, 876, 1928, 1069,\\
	2321, 1269, 2763, 1513, 3259, 1765, 3810, 2066, 4421, 2376, 5093, 2740, \dots 
\end{eqnarray*}

This is sequence A122133 in the \textit{Online-Encyclopedia of Integer Sequences} \cite{sloane}.
	
\medskip

Besides the maximum area $n+\left\lfloor\frac{n-1}{2}\right\rfloor\left\lfloor\frac{n}{2}\right\rfloor$ and 
the minimum area $n$ of the convex hull of polyominoes with $n$ cells several other values may be attained. In 
Theorem \ref{thm_possible_values} we have completely characterized the set of areas of the convex hull of polyominoes 
with $n$ cells.

\medskip

{\sc Proof of Theorem \ref{thm_possible_values}.} 
	Since the vertex points of the convex hull of a polyomino are lattice points on an integer grid the area of the convex 
	hull is an integral multiple of $\frac{1}{2}$. with Theorem \ref{thm_maximum_2} we conclude that the desired set is a 
	subset of 
	$$
		S=\left\{n+\frac{m}{2}\mid m\le\left\lfloor\frac{n-1}{2}\right\rfloor\left\lfloor\frac{n}{2}\right\rfloor,
		\,m\in\mathbb{N}_0\right\}.
	$$
	A polyomino $P$ consisting of $n$ cells with area $n+\frac{1}{2}$ of the convex hull must contain a triangle of area 
	$\frac{1}{2}$. If we extend the triangle to a square we get a convex polyomino $P'$ consisting of $n+1$ cells. Thus 
	$P'$ is an rectangular $s\times t$-polyomino with $s\cdot t=n+1$ and $s,t\in\mathbb{N}$. If $n+1$ is a prime there 
	exists only the $1\times (n+1)$-polyomino where deleting a square yields an area of $n$ for the convex hull. So we 
	have to exclude this case in the above set $S$ and receive the proposed set.
	
	\begin{figure}[h]
		\begin{center}
			\setlength{\unitlength}{0.5cm}
			\begin{picture}(8,4)
				\put(4,0){$a$}
				\put(0,0.6){\line(1,0){8}}
				\put(0,0.6){\line(0,1){0.2}}
				\put(8,0.6){\line(0,1){0.2}}
				\put(0,1){\line(1,0){8}}
				\put(0,2){\line(1,0){8}}
				\put(0,3){\line(1,0){5}}
				\put(6,3){\line(1,0){1}}
				\put(0,1){\line(0,1){2}}
				\put(1,1){\line(0,1){2}}
				\put(2,1){\line(0,1){2}}
				\put(3,1){\line(0,1){2}}
				\put(4,1){\line(0,1){2}}
				\put(5,1){\line(0,1){2}}
				\put(6,1){\line(0,1){2}}
				\put(7,1){\line(0,1){2}}
				\put(8,1){\line(0,1){1}}
				\put(0,3.4){\line(0,-1){0.2}}
				\put(5,3.4){\line(0,-1){0.2}}
				\put(6,3.4){\line(0,-1){0.2}}
				\put(7,3.4){\line(0,-1){0.2}}
				\put(0,3.4){\line(1,0){5}}
				\put(6,3.4){\line(1,0){1}}
				\put(6.15,2.2){$\leftrightarrow$}
				\put(2.4,3.8){$b$}
				\put(6.35,3.8){$1$}
				\put(5.2,2.4){\line(1,0){0.6}}
				\put(5.4,2.6){$l$}
			\end{picture}
			\caption{Construction 1: $2\le m\le 2n-8$.}
			\label{fig_construction_1}
		\end{center}
	\end{figure}
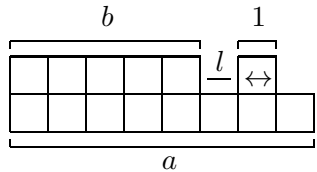
		
	For the other direction we give some constructions. For $m=0$ we have the rectangular $1\times n$-polyomino as an 
	example. The above consideration $m=1$ yields  a construction if $n+1$ is a composite number. Now we consider 
	Construction 1 depicted in Figure \ref{fig_construction_1}. We choose $n=a+b+1$, 
	$\left\lceil\frac{n}{2}\right\rceil\le a\le n-2$, and $0\le l\le a-b-1=2a-n$. Thus $a\ge b+1$ and Construction 1 
	is possible. If we run through the possible values of $a$ and $l$ we obtain examples for
	\begin{eqnarray*} 
		m\in&&\{0\},\{2,3,4\},\dots,\{2a-n,\dots,4a-2n\},\dots,\{n-4,\dots,2n-8\}\\
		&=&\{0,2,3,\dots,2n-8\}
	\end{eqnarray*}
	if $n\equiv 0\mod 2$ and for
	\begin{eqnarray*}
		m\in&&\{1,2\},\{3,4,5,6\},\dots,\{2a-n,\dots,4a-2n\},\dots,\{n-4,\dots,2n-8\}\\
		&=& \{1,2,\dots,2n-8\}
	\end{eqnarray*}
	if $n\equiv 1\mod 2$.
	
	\begin{figure}[h]
		\begin{center}
			\setlength{\unitlength}{0.5cm}
			\begin{picture}(13,5.5)
				\put(1,0){\line(1,0){4}}
				\put(0,1){\line(1,0){5}}
				\put(0,2){\line(1,0){2}}
				\put(0,3){\line(1,0){1}}
				\put(0,1){\line(0,1){2}}
				\put(1,0){\line(0,1){3}}
				\put(2,0){\line(0,1){2}}
				\put(3,0){\line(0,1){1}}
				\put(4,0){\line(0,1){1}}
				\put(5,0){\line(0,1){1}}
				\put(9,0){\line(0,1){2}}
				\put(10,0){\line(0,1){2}}
				\put(11,0){\line(0,1){1}}
				\put(12,0){\line(0,1){1}}
				\put(13,0){\line(0,1){1}}
				\put(9,0){\line(1,0){4}}
				\put(6,1){\line(1,0){7}}
				\put(6,2){\line(1,0){4}}
				\put(6,3){\line(1,0){1}}
				\put(6,1){\line(0,1){2}}
				\put(7,1){\line(0,1){2}}
				\put(8,1){\line(0,1){1}}
				\put(9,1){\line(0,1){1}}
				\put(1,3.5){\line(0,1){1}}
				\put(2,3.5){\line(0,1){1}}
				\put(3,2.5){\line(0,1){2}}
				\put(4,2.5){\line(0,1){2}}
				\put(5,2.5){\line(0,1){1}}
				\put(1,3.5){\line(1,0){4}}
				\put(1,4.5){\line(1,0){3}}
				\put(3,2.5){\line(1,0){2}}
				\put(7,3.5){\line(0,1){2}}
				\put(8,2.5){\line(0,1){3}}
				\put(9,2.5){\line(0,1){2}}
				\put(10,2.5){\line(0,1){1}}
				\put(11,2.5){\line(0,1){1}}
				\put(7,3.5){\line(1,0){4}}
				\put(7,4.5){\line(1,0){2}}
				\put(7,5.5){\line(1,0){1}}
				\put(8,2.5){\line(1,0){3}}
			\end{picture}
			\quad\,\,\,\,
			\begin{picture}(9,3)
				\put(0,0){\line(1,0){1}}
				\put(0,1){\line(1,0){6}}
				\put(0,2){\line(1,0){6}}
				\put(0,3){\line(1,0){1}}
				\put(3,3){\line(1,0){1}}
				\put(4,0){\line(1,0){1}}
				\put(0,0){\line(0,1){1}}
				\put(1,0){\line(0,1){1}}
				\put(4,0){\line(0,1){1}}
				\put(5,0){\line(0,1){1}}
				\put(0,1){\line(0,1){1}}
				\put(1,1){\line(0,1){1}}
				\put(2,1){\line(0,1){1}}
				\put(3,1){\line(0,1){1}}
				\put(4,1){\line(0,1){1}}
				\put(5,1){\line(0,1){1}}
				\put(6,1){\line(0,1){1}}
				\put(0,2){\line(0,1){1}}
				\put(1,2){\line(0,1){1}}
				\put(3,2){\line(0,1){1}}
				\put(4,2){\line(0,1){1}}
				\put(6.4,1.4){$\dots$}
				\put(7.8,1){\line(1,0){1}}
				\put(7.8,2){\line(1,0){1}}
				\put(8.8,1){\line(0,1){1}}
				\put(7.8,1){\line(0,1){1}}
			\end{picture}
			\caption{Construction 2: $m=2n-7$.}
			\label{fig_construction_2}
		\end{center}
	\end{figure}
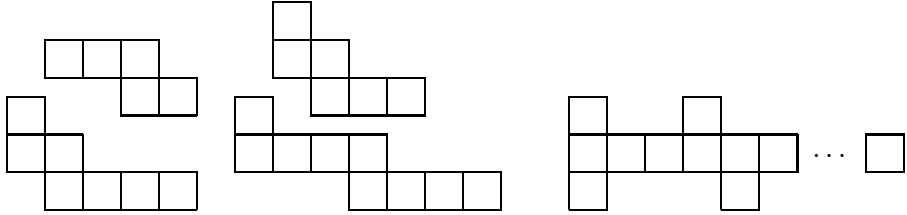
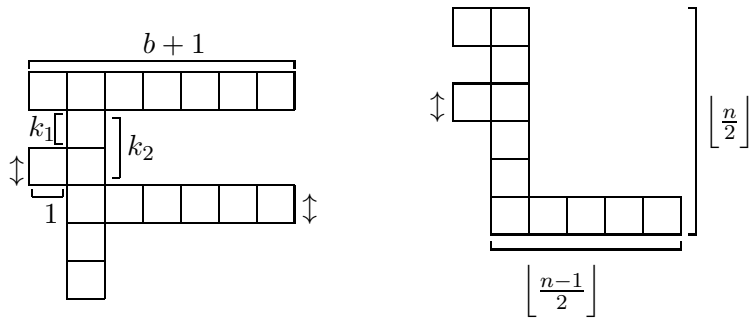
\begin{figure}[h]
		\begin{center}
			\setlength{\unitlength}{0.5cm}
			\begin{picture}(9,6)
				\put(1,6){\line(1,0){7}}
				\put(1,5){\line(1,0){7}}
				\put(2,4){\line(1,0){1}}
				\put(2,3){\line(1,0){1}}
				\put(2,2){\line(1,0){1}}
				\put(2,1){\line(1,0){1}}
				\put(2,0){\line(1,0){1}}
				\put(2,6){\line(0,-1){6}}
				\put(3,6){\line(0,-1){6}}
				\put(1,6){\line(0,-1){1}}
				\put(3,6){\line(0,-1){1}}
				\put(4,6){\line(0,-1){1}}
				\put(5,6){\line(0,-1){1}}
				\put(6,6){\line(0,-1){1}}
				\put(7,6){\line(0,-1){1}}
				\put(8,6){\line(0,-1){1}}
				\put(2,4){\line(-1,0){1}}
				\put(2,3){\line(-1,0){1}}
				\put(1,4){\line(0,-1){1}}
				\put(3,3){\line(1,0){5}}
				\put(3,2){\line(1,0){5}}
				\put(4,3){\line(0,-1){1}}
				\put(5,3){\line(0,-1){1}}
				\put(6,3){\line(0,-1){1}}
				\put(7,3){\line(0,-1){1}}
				\put(8,3){\line(0,-1){1}}
				\put(3.4,4.8){\line(0,-1){1.6}}
				\put(3.4,4.8){\line(-1,0){0.2}}
				\put(3.4,3.2){\line(-1,0){0.2}}
				\put(3.6,3.8){$k_2$}
				\put(1.0,4.3){$k_1$}
				\put(1.7,4.1){\line(0,1){.8}}
				\put(1.7,4.1){\line(1,0){.2}}
				\put(1.7,4.9){\line(1,0){.2}}
				\put(8.2,2.2){$\updownarrow$}
				\put(.5,3.2){$\updownarrow$}
				\put(1,6.3){\line(1,0){7}}
				\put(1,6.3){\line(0,-1){0.2}}
				\put(8,6.3){\line(0,-1){0.2}}
				\put(4,6.5){$b+1$}
				\put(1.1,2.7){\line(1,0){.8}}
				\put(1.1,2.7){\line(0,1){0.2}}
				\put(1.9,2.7){\line(0,1){0.2}}
				\put(1.4,2){$1$}
			\end{picture}
  			\quad\quad\quad
			\begin{picture}(8.4,7.7)
				\put(1.6,1.7){\line(1,0){5}}
				\put(1.6,2.7){\line(1,0){5}}
				\put(1.6,3.7){\line(1,0){1}}
				\put(0.6,4.7){\line(1,0){2}}
				\put(0.6,5.7){\line(1,0){2}}
				\put(0.6,6.7){\line(1,0){2}}
				\put(0.6,7.7){\line(1,0){2}}
				\put(1.6,1.7){\line(0,1){6}}
				\put(2.6,1.7){\line(0,1){6}}
				\put(0.6,6.7){\line(0,1){1}}
				\put(0.6,4.7){\line(0,1){1}}
				\put(0,4.9){$\updownarrow$}
				\put(3.6,1.7){\line(0,1){1}}
				\put(4.6,1.7){\line(0,1){1}}
				\put(5.6,1.7){\line(0,1){1}}
				\put(6.6,1.7){\line(0,1){1}}
				\put(1.6,1.3){\line(1,0){5}}
				\put(1.6,1.3){\line(0,1){0.2}}
				\put(6.6,1.3){\line(0,1){0.2}}
				\put(2.4,0){$\Big\lfloor\frac{n-1}{2}\Big\rfloor$}
				\put(7.0,1.7){\line(0,1){6}}
				\put(7.0,1.7){\line(-1,0){0.2}}
				\put(7.0,7.7){\line(-1,0){0.2}}
				\put(7.2,4.5){$\Big\lfloor\frac{n}{2}\Big\rfloor$}
			\end{picture}
			\caption{Construction 3 and Construction 4.}
			\label{fig_construction_3_and_4}
		\end{center}
	\end{figure}

	In Figure \ref{fig_construction_2} we give a construction for $m=2n-7$ and in Figure \ref{fig_construction_3_and_4} we give on the 
	left hand side a construction for $2n-6\le m\le \left\lceil\frac{n^2-4n}{4}\right\rceil$ with parameters $k_1$, $k_2$, and $b$. 
	The conditions for these parameters are $0\le k_1,k_2\le n-2b-2$ and $n-2b-2\ge b$. With given $k_1,k_2,b,n$ we have 
	$m=bn-2b^2-2b+k_1+k_2(b-1)$. 
	Since we can vary $k_1$ at least between $0$ and $b-1$ we can produce for a fix $b$ all values 
	beetween $b(n-2b-2)$ and $2b(n-2b-2)$ by varying $k_1$ and $k_2$.  Now we want to combine those 
	intervals for successive values for $b$. The assumption that the intervals leave a gap is 
	equivalent to $2(b-1)(n-2(b-1)-2)<b(n-2b-2)$, that is, $n < 2b\frac{b-3}{b-2}$. 
	We choose $2\le b\le \left\lfloor\frac{n}{4}\right\rfloor$ and receive constructions for 
	$$
		m\in\left\{2n-6,2n-5,\dots,\left\lceil\frac{n^2-4n}{4}\right\rceil\right\}.
	$$
	
	\bigskip
	
	On the right hand side of Figure \ref{fig_construction_3_and_4} we give a construction for $n\ge 5$ and 
	$$
		m\in\left\{\left\lfloor\frac{n^2-4n}{4}\right\rfloor,\dots,\left\lfloor\frac{n^2-2n-8}{4}\right\rfloor\right\}.
	$$
	
	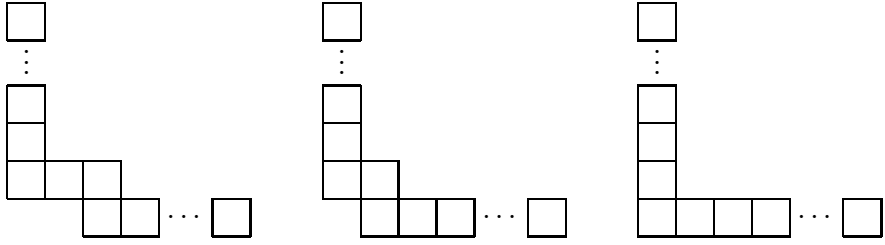
\begin{figure}[h]
		\begin{center}
			\setlength{\unitlength}{0.5cm}
			\begin{picture}(6.5,6.3)
				\put(0,1){\line(0,1){3}}
				\put(1,1){\line(0,1){3}}
				\put(2,0){\line(0,1){2}}
				\put(3,0){\line(0,1){2}}
				\put(4,0){\line(0,1){1}}
				\put(4.2,0.5){$\dots$}
				\put(0,4){\line(1,0){1}}
				\put(0,3){\line(1,0){1}}
				\put(0,2){\line(1,0){3}}
				\put(0,1){\line(1,0){4}}
				\put(2,0){\line(1,0){2}}
				\put(5.4,0){\line(0,1){1}}
				\put(6.4,0){\line(0,1){1}}
				\put(5.4,0){\line(1,0){1}}
				\put(5.4,1){\line(1,0){1}}
				\put(0.4,4.3){$\vdots$}
				\put(0,5.2){\line(1,0){1}}
				\put(0,6.2){\line(1,0){1}}
				\put(0,5.2){\line(0,1){1}}
				\put(1,5.2){\line(0,1){1}}
			\end{picture}
			\quad\quad
			\begin{picture}(6.5,6.3)
				\put(0,1){\line(0,1){3}}
				\put(1,1){\line(0,1){3}}
				\put(2,0){\line(0,1){2}}
				\put(3,0){\line(0,1){1}}
				\put(1,0){\line(0,1){1}}
				\put(4,0){\line(0,1){1}}
				\put(4.2,0.5){$\dots$}
				\put(0,4){\line(1,0){1}}
				\put(0,3){\line(1,0){1}}
				\put(0,2){\line(1,0){2}}
				\put(0,1){\line(1,0){4}}
				\put(1,0){\line(1,0){3}}
				\put(5.4,0){\line(0,1){1}}
				\put(6.4,0){\line(0,1){1}}
				\put(5.4,0){\line(1,0){1}}
				\put(5.4,1){\line(1,0){1}}
				\put(0.4,4.3){$\vdots$}
				\put(0,5.2){\line(1,0){1}}
				\put(0,6.2){\line(1,0){1}}
				\put(0,5.2){\line(0,1){1}}
				\put(1,5.2){\line(0,1){1}}
			\end{picture}
			\quad\quad
			\begin{picture}(6.5,6.3)
				\put(0,0){\line(0,1){4}}
				\put(1,1){\line(0,1){3}}
				\put(2,0){\line(0,1){1}}
				\put(3,0){\line(0,1){1}}
				\put(1,0){\line(0,1){1}}
				\put(4,0){\line(0,1){1}}
				\put(4.2,0.5){$\dots$}
				\put(0,4){\line(1,0){1}}
				\put(0,3){\line(1,0){1}}
				\put(0,2){\line(1,0){1}}
				\put(0,1){\line(1,0){4}}
				\put(0,0){\line(1,0){4}}
				\put(5.4,0){\line(0,1){1}}
				\put(6.4,0){\line(0,1){1}}
				\put(5.4,0){\line(1,0){1}}
				\put(5.4,1){\line(1,0){1}}
				\put(0.4,4.3){$\vdots$}
				\put(0,5.2){\line(1,0){1}}
				\put(0,6.2){\line(1,0){1}}
				\put(0,5.2){\line(0,1){1}}
				\put(1,5.2){\line(0,1){1}}
			\end{picture}
			\caption{Construction 5: $\left\lfloor\frac{n^2-2n-6}{4}\right\rfloor\le m\le \left\lfloor\frac{n^2-2n+2}{4}\right\rfloor$.}
			\label{fig_construction_5}
		\end{center}
		Constructions for the remaining values 
		$$
			m\in\left\{
			\left\lfloor\frac{n^2-2n-6}{4}\right\rfloor,
			\left\lfloor\frac{n^2-2n-2}{4}\right\rfloor,
			\left\lfloor\frac{n^2-2n+2}{4}\right\rfloor=
			\left\lfloor\frac{n-1}{2}\right\rfloor\left\lfloor\frac{n}{2}\right\rfloor
			\right\}
		$$
		are given in Figure \ref{fig_construction_5}.
	\end{figure}
	
\hfill{$\square$}

\section{Dimensions $\mathbf{d\ge 3}$}
\label{sec_3}

To prove the conjecture of Bezdek, Bra\ss, and Harborth for dimensions $d\ge 3$ we proceed similar as in Section \ref{sec_2}. 

\begin{definition}
	\label{main_definition}
	$$
		f_d(l_1,\dots ,l_d,v_1,\dots, v_d)=\sum_{I\subseteq\{1,\dots ,d\}}\frac{1}{|I|!2^{d-|I|}}
		\sum_{b=0}^{2^d-1}\prod_{i\in I}q_{b,i}
	$$
	with $d\ge 1$ and $b=\sum\limits_{j=1}^db_j2^{j-1}$, $b_j\in\{0,1\}$,
	$q_{b,i}=\left\{\begin{array}{ccc}l_i-1 &for& b_i=0\, , \\v_i &for& b_i=1\, . \end{array}\right.$
\end{definition}

\begin{example}
	\begin{eqnarray*}
		\!\!\!\!\!\!\!\!&&f_3(l_1,l_2,l_3,v_1,v_2,v_3)=1+(l_1-1)+(l_2-1)+(l_3-1)+\frac{(l_1-1)(l_2-1)}{2}+\\
      \!\!\!\!\!\!\!\!&&\frac{(l_1-1)(l_3-1)}{2}+\frac{(l_2-1)(l_3-1)}{2}+\frac{(l_1-1)(l_2-1)(l_3-1)}{6}+\frac{v_1(l_2-1)}{2}+\\
      \!\!\!\!\!\!\!\!&&\frac{v_1(l_3-1)}{2}+\frac{v_2(l_1-1)}{2}+\frac{v_2(l_3-1)}{2}+\frac{v_3(l_1-1)}{2}+
		\frac{v_3(l_2-1)}{2}+\\
      \!\!\!\!\!\!\!\!&&\frac{v_1(l_2-1)(l_3-1)}{6}+\frac{v_2(l_1-1)(l_3-1)}{6}+\frac{v_3(l_1-1)(l_2-1)}{6}+\frac{v_1v_2(l_3-1)}{6}+\\
      \!\!\!\!\!\!\!\!&&\frac{v_1v_3(l_2-1)}{6}+\frac{v_2v_3(l_1-1)}{6}+v_1+v_2+v_3+\frac{v_1v_2}{2}+\frac{v_1v_3}{2}+\frac{v_2v_3}{2}
		+\frac{v_1v_2v_3}{6}\, .
	\end{eqnarray*}
\end{example}

\begin{lemma}
	\label{main_lemma}
	The $d$-dimensional volume of the convex hull of a polyomino with $n$ unit hypercubes is at most 
	$f_d(l_1,\dots,l_d,v_1,\dots,v_d)$.
\end{lemma}
\begin{proof}
	We prove the statement by double induction on $d$ and $n$, using equation \ref{main_equation}. Since 
	the case $d=2$ is already done in Theorem \ref{thm_maximum_2} we assume that the lemma is proven for 
	the $\overline{d}<d$. Since for $n=1$ only $l_i=1$, $v_i=0$, $i\in\{1,\dots,d\}$ is possible and 
	$f_d(1,\dots,1,0,\dots,0)=1$ the induction base for $n$ is done. Now we assume that the lemma is proven 
	for all possible tuples $(l_1,\dots,l_d,v_1,\dots,v_d)$ with $1+\sum_{i=1}^d(l_i-1)+\sum_{i=1}^d v_i=n-1$. 
	Due to symmetry we consider only the growth of $l_1$ or $v_1$, and the volume of the convex hull by 
	adding the $n$-th hypercube.
	\begin{figure}[h]
		\begin{center}
			\setlength{\unitlength}{0.36cm}
			\begin{picture}(34,10)
				\put(5,3){\line(1,0){4}}
				\put(3,5){\line(0,1){4}}
				\put(3,5){\line(1,0){6}}
				\put(5,3){\line(0,1){6}}
				\put(9,3){\line(0,1){2}}
				\put(3,9){\line(1,0){2}}
				\qbezier[20](5,9)(9,9)(9,5)
				\qbezier[15](4,4)(4.2,6)(8,8)
				\put(2,4){\line(1,5){1}}
				\put(4,4){\line(1,5){1}}
				\put(4,2){\line(5,1){5}}
				\put(4,4){\line(5,1){5}}
				\thicklines
				\put(2.0,2.0){\line(1,0){2}}
				\put(2.0,2.0){\line(0,1){2}}
				\put(4.0,4.0){\line(-1,0){2}}
				\put(4.0,4.0){\line(0,-1){2}}
				\put(2,4){\line(1,1){1}}
				\put(4,4){\line(1,1){1}}
				\put(4,2){\line(1,1){1}}
				\put(2,2){\line(1,1){2}}
				\put(3,2){\line(1,1){1}}
				\put(2,3){\line(1,1){1}}
				\put(5,0.8){$l_3\,\rightarrow$}
				\put(0.4,6){$l_2\,\uparrow$}
				\put(7,6){$A$}
				\thinlines
				\put(13,5){\line(1,0){3}}
				\put(13,7){\line(1,0){8}}
				\put(13,5){\line(0,1){2}}
				\put(21,5){\line(0,1){2}}
				\put(17,7){\line(0,1){2}}
				\put(19,7){\line(0,1){2}}
				\put(17,9){\line(1,0){2}}
				\put(17,7){\line(-1,-1){1}}
				\put(19,7){\line(-1,-1){1}}
				\put(18,4){\line(1,1){1}}
				\put(19,5){\line(0,1){2}}
				\put(19,5){\line(1,0){2}}
				\put(19,5){\line(0,-1){4}}
				\put(19,1){\line(-1,0){2}}
				\put(17,1){\line(0,1){3}}
				\put(16,6){\line(-3,1){3}}
				\put(16,4){\line(-3,1){3}}
				\put(18,6){\line(3,1){3}}
				\put(18,4){\line(3,1){3}}
				\put(16,6){\line(1,3){1}}
				\put(18,6){\line(1,3){1}}
				\put(16,4){\line(1,-3){1}}
				\put(18,4){\line(1,-3){1}}
				\qbezier[18](13,7)(13,9)(17,9)
				\qbezier[15](19,9)(21,9)(21,7)
				\qbezier[20](13,5)(13,1)(17,1)
				\qbezier[18](19,1)(21,1)(21,5)
				\thicklines
				\put(16,6){\line(1,0){2}}
				\put(16,6){\line(0,-1){2}}
				\put(18,4){\line(-1,0){2}}
				\put(18,4){\line(0,1){2}}
				\put(16,4){\line(1,1){2}}
				\put(16,5){\line(1,1){1}}
				\put(17,4){\line(1,1){1}}
				\thinlines
				\put(30,7){\line(1,0){2}}
				\put(30,5){\line(1,0){2}}
				\put(32,5){\line(0,1){2}}
				\put(30,5){\line(0,1){2}}
				\put(27,2){\line(5,3){5}}
				\put(27,4){\line(5,3){5}}
				\qbezier[5](28,7)(29,7)(30,7)
				\qbezier[15](28,3)(32,3)(32,5)
				\put(26,7){\line(1,0){2}}
				\put(26,7){\line(-1,-3){1}}
				\put(28,7){\line(-1,-3){1}}
				\put(26,5){\line(0,1){2}}
				\put(28,5){\line(0,-1){2}}
				\put(28,5){\line(0,1){2}}
				\put(26,5){\line(1,0){2}}
				\put(25,4){\line(1,1){1}}
				\put(27,4){\line(1,1){1}}
				\put(27,2){\line(1,1){1}}
				\thicklines
				\put(25.0,2.0){\line(1,0){2}}
				\put(25.0,2.0){\line(0,1){2}}
				\put(27.0,4.0){\line(-1,0){2}}
				\put(27.0,4.0){\line(0,-1){2}}
				\put(25,2){\line(1,1){2}}
				\put(26,2){\line(1,1){1}}
				\put(25,3){\line(1,1){1}}
				\put(28,0.8){$l_3\,\rightarrow$}
				\put(23.4,6){$l_2\,\uparrow$}
			\end{picture}
			\caption{Increasing $l_1$ in the 3-dimensional case.}
			\label{fig_increasing_3}
		\end{center}
	\end{figure}
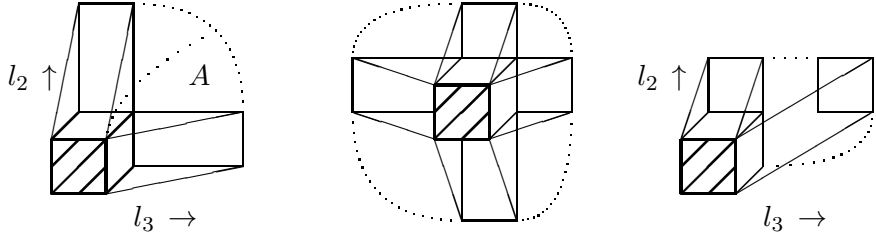
	
	\vspace*{-8mm}
	
	\begin{enumerate}
		\item[(i)]	$l_1$ increases by one:
		
						As in the proof of Lemma \ref{lemma_two_dimensional} draw lines of the convex hull of the $n$-th cube 
						and its neighbor cube $N$, see Figure \ref{fig_increasing_3} for a $3$-dimensional example. To be more 
						precisely each line of the new convex hull has a corner point $X$ of the upper face of the $n$-th cube 
						as an endpoint. We will denote the second endpoint of this line by $Y$. In direction of axis $1$ there 
						is a corner point $\overline{X}$ of the bottom face of the $n$-th cube. Since $\overline{X}$ is also a 
						corner point of $N$ the line $\overline{X}Y$ is part of the old convex hull if $Y$ is part of the old 
						convex hull. In this case we draw the line $\overline{X}Y$. In the other case $Y$ is also a corner point 
						of the upper face of the new cube and we draw the line $\overline{X}\overline{Y}$ where $\overline{Y}$ 
						is similar defined as $\overline{X}$. Additionally we draw all lines $XY$ and $X\overline{X}$. 
						
						Doing this we have constructed a geometrical body which contains the increase of the convex hull and is 
						subdivided into nice geometrical objects $O_i$ with volume $\frac{\text{base}\,\times\,\text{height}}{k_i}$, 
						for some $k_i\in\{1,\dots,d\}$ each. For dimension $d=3$ the cases $k_i=1$, $k_i=2$, or $k_i=3$ correspond 
						to a box, a prism, or a tetrahedron.
						
						We project the convex hull of the whole polyomino into the hyperplane orthogonal to axis direction $1$ and 
						receive a hypervolume $A$. This is the convex hull of a $(d-1)$-dimensional polyomino with parameters 
						$\overline{l}_2,\dots,\overline{l}_d,\overline{v}_2,\dots,\overline{v}_d$ where $\overline{l}_i\le l_i$ and 
						$\overline{v}_i\le v_i$. From the induction hypothesis we know $A\le f_{d-1}(l_2,\dots,l_d,v_2,\dots,v_d)$. 
						We apply the same projection to the $O_i$ and objects $A_i$. Due to the construction the $A_i$ are non 
						overlapping and we have $\sum A_i\le A$. Using Cavalieri's theorem we determine the volume of $O_i$ to be 
						$\frac{A_i\times 1}{k_i}$. More precisely, we choose lines of the form $\overline{X}X$ as height and lift 
						the old base up until it is orthogonal to axis direction $1$. Thus we may assign a factor $\frac{1}{k}$ to 
						each piece of $A$ to bound the growth of the volume of the convex hull. We estimate the parts in a way that 
						the parts with the higher factors are as big as theoretical possible. 
						
						For every $0\le r \le d-1$ we consider the sets $\{i_1,i_2,\dots ,i_r\}$ with $1\ne i_a\ne i_b$ for $a\ne b$. 
						Let $Z$ be such a set. Define $\overline{Z}=\{j_1,\dots,j_{d-r-1}\}$ by $Z\cap\overline{Z}=\{\}$ and $Z\cup 
						\overline{Z}=\{2,\dots ,d\}$. So the vector space spanned by the axis directions of $Z$ and the vector space 
						spanned by the axis directions of $\overline{Z}$ are orthogonal. If we project the convex hull in the vector 
						space spanned by $\overline{Z}$ the resulting volume is at most$f_{d-r-1}(l_{j_1},\dots,l_{j_{d-r-1}},
						v_{j_1},$ $\dots,v_{j_{d-r-1}})$ since it is the convex hull of a $(d-r-1)$-dimensional polyomino. Since 
						$\overline{Z}$ has cardinality $d-r-1$ the set $Z$ yields a contribution of $\frac{1}{d-r}f_{d-r-1}(l_{j_1},
						\dots,l_{j_{d-r-1}},v_{j_1},\dots,v_{j_{d-r-1}})$ to the volume of the convex hull. With the notations 
						from Definition \ref{main_definition} this is
						$$
							\frac{1}{d-r}\sum_{I\subseteq\{j_1,\dots ,j_{r-d-1}\}}\frac{1}{|I|!2^{d-r-1-|I|}}
							\sum_{b=0}^{2^{d-r-1}-1}\prod_{i\in I}q_{b,i}\, .
						$$
						Our aim is to assign the maximum possible factor to each part of  $A$. For that reason we count for $Z$ a 
						maximum contribution of 
						$$
							\frac{1}{d-r}\frac{1}{|d-r-1|!}
							\sum_{b=0}^{2^{d-r-1}-1}\prod_{i\in \overline{Y}}q_{b,i}
						$$
						to the volume of the convex hull.
						
						If we do so for all possible sets $Z$ we have assigned a factor between $1$ and $\frac{1}{d}$ to every 
						summand of $f_{d-1}(l_2,\dots ,l_d,v_2,\dots,v_d)$. To get the induction step now we have to remark that 
						the above described  sum with its factors is exactly the difference between  $f_{d}(l_1+1,\dots ,l_d,v_1,
						\dots,v_d)$ and $f_{d}(l_1,\dots ,l_d,v_1,\dots,v_d)$.
		\item[(ii)]	$v_1$ increases by one:
						
						Due to symmetry of the $l_i$ and $v_i$ in Definition \ref{main_definition} this is similar to case (i). 
						Additionally we remark that the maximum cannot be achieved in this case since we double count a part 
						of the contribution of the new cube to the volume of the convex hull in our estimations.
	\end{enumerate}
\end{proof}

\begin{theorem}
	\label{main_theorem}
	The $d$-dimensional volume of the convex hull of any facet-to-facet connected system of $n$ unit
	hypercubes is 
	$$
		\sum\limits_{I\subseteq\{1,\dots,d\}}\frac{1}{|I|!}\prod\limits_{i\in
		I}\left\lfloor\frac{n-2+i}{d}\right\rfloor\, .
	$$
\end{theorem}
\begin{proof}
	For given $n$ we determine the maximum of  $f_d(l_1,\dots,f_d,v_1,\dots,v_d)$. 
	Due to $$f_d(l_1+1,l_2,\dots ,l_d,v_1-1,v_2,\dots ,v_d)-f_d(l_1,l_2,\dots ,l_d,v_1,v_2,\dots
	,v_d)=0$$ and due to symmetry we assume $v_1=\dots=v_d=0$ and $l_1\le l_2\le\dots\le l_d$. Since
	\begin{equation}
		\label{ieq_main_thm}
		f_d(l_1+1,l_2,\dots ,l_{d-1},l_d-1,0,0,\dots ,0)-f_d(l_1,l_2,\dots ,l_d,0,0,\dots,0)>0
	\end{equation}
	we have $0\le l_d-l_1\le 1$. Inequality \ref{ieq_main_thm} due to the following consideration.
	If a summand of $f_d(\dots)$ contains the term $l_1$ and does not contain $l_d$ then there will 
	be a corresponding summand with $l_1$ replaced by $l_d$, so those terms equalize each other in the 
	above difference. Clearly the summands containing none of the terms $l_1$ or $l_d$ equalize each other 
	in the difference. So there are left only the summands with both terms $l_1$ and $l_d$. Since
	$(l_1+1-1)(l_d-1-1)-(l_1-1)(l_d-1)=l_d-l_1-1>0$ inequality \ref{ieq_main_thm} is valid. 
	
	\medskip
	
	Combining equation \ref{main_equation} with $0\le l_d-l_1\le 1$ and $l_1\le l_2\le\dots\le l_d$ gives 
	$l_i=\left\lfloor\frac{n-2+i+d}{d}\right\rfloor$. Thus by inserting in Lemma \ref{main_lemma} we receive 
	the upper bound. The maximum is attained for example by a polyomino consisting of $d$ pairwise orthogonal
	linear arms with $\left\lfloor\frac{n-2+i}{d}\right\rfloor$ cubes ($i=1\dots d$) joined to a central cube.
\end{proof}

\begin{conjecture}
	Every $d$-dimensional polyomino $P$ with parameters $l_1,\dots,l_d,$ $v_1,\dots,v_d$ and maximum volume 
	of the convex hull fulfills $v_1=\dots=v_d=0$ and contains a sub polyomino $P'$ fulfilling:
	\begin{enumerate}
		\item[(i)]	$P'$ has height $1$ in direction of axis $i$,
		\item[(ii)]	the projection of $P'$ along $i$ has also maximal volume of the convex hull and 
						parameters $l_1,\dots,l_{i-1},l_{i+1},\dots,l_d$,
		\item[(iii)] $P$ can be decomposed into $P'$ and up to two orthogonal linear arms.
	\end{enumerate}
\end{conjecture}

We remark that $v_1=\dots=v_d=0$ and the maximality of the volume of the convex hull of sub polyominoes and projections 
of $P$ can be concluded from the proof of Theorem \ref{main_theorem}.

\begin{lemma}
	If there exists a $d$-dimensional polyomino with $n$ cells and volume $v$ of the convex hull, then  $v\in V_{d,n}$ with
	$$
		V_{d,n}=\left\{n+\frac{m}{d!}\Big|m\le \sum\limits_{I\subseteq\{1,\dots,d\}}\frac{d!}{|I|!}\prod\limits_{i\in
		I}\left\lfloor\frac{n-2+i}{d}\right\rfloor\,m\in \mathbb{N}_0\right\}\, .
	$$
\end{lemma}  
\begin{proof}
	For the determination of the volume of the convex hull of a $d$-dimensional polyomino we only have to consider the 
	set of $S$ corner points of its hypercubes which lie on an integer grid. We can decompose the convex hull into 
	$d$-dimensional simplices with the volume 
	$$
		\frac{1}{d!}\left|\begin{array}{cccc}    
		x_{1,1} & \dots & x_{1,d} & 1 \\
		\vdots & \ddots & \vdots & \vdots \\
		x_{d+1,1} & \dots & x_{d+1,d} & 1 
		\end{array}\right|
	$$
	where the coordinates of the $d+1$ points are given by $(x_{i,1},\dots,x_{i,d})\in\mathbb{Z}^d$. Thus the volume of 
	the convex hull is an integer multiple of $\frac{1}{d}$. The lower bound $n\le v$ is obvious and the upper bound is 
	given by Theorem \ref{main_theorem}.
\end{proof}

\section{Remarks}
\label{sec_remark} 

We leave the description and the enumeration of the polyominoes with maximum convex hull for dimension $d\ge 3$ as a 
task for the interested reader. It would also be nice to see a version of Theorem \ref{thm_possible_values} for 
higher dimensions.

\medskip

The authors of \cite{Bezdek1994} mention another class of problems which are related to the problems in \cite{circles} and \cite{huellen}: What is the maximum area of the convex hull of all connected edge-to-edge packings of $n$ congruent regular 
$k$-gons (also denoted as $k$-polyominoes, see \cite{kpolyominoes}) in the plane. The methods of Section \ref{sec_2} might 
be applicable for these problems. 

\begin{conjecture}
	The area of the convex hull of any edge-to-edge connected system
	of regular unit hexagons is at most $\frac{1}{6}\left\lfloor n^2+\frac{14}{3}n+1\right\rfloor$.
\end{conjecture}

\nocite{kurz}
\bibliography{convex_hulls}
\bibdata{convex_hulls}
\bibliographystyle{plain}  

\end{document}